\documentclass[12pt]{article}
\usepackage{amsfonts}
\usepackage{amssymb,amsfonts,amsmath,amsthm,cite,color}
\usepackage{enumerate}
\usepackage{array}

\newcommand{\supp}{\hbox{\rm supp}}

\PassOptionsToPackage{hyphens}{url}
\usepackage{hyperref}
\allowdisplaybreaks[4]

\newtheorem{thm}{Theorem}[section]

\newtheorem{lem}[thm]{Lemma}
\newtheorem{pro}[thm]{Proposition}
\newtheorem{con}[thm]{Conjecture}
\newtheorem{exa}[thm]{Example}

\makeatletter \@addtoreset{equation}{section}

\def\pf{\noindent {\it Proof.\ }}
\def\qed{\hfill \rule{4pt}{7pt}}

\title{\vspace{-1.6cm}
Powerful sets: a generalisation of binary matroids\thanks{This work was presented at the 40th Australasian Conference on Combinatorial Mathematics and Combinatorial Computing (40ACCMCC), University of Newcastle, Australia, Dec.\ 2016.}}
\author{Graham E. Farr   \\
Faculty of Information Technology   \\
Monash University   \\
Clayton, Victoria 3800   \\
Australia   \\
Email: \href{mailto:Graham.Farr@monash.edu}{\texttt{Graham.Farr@monash.edu}}
\and
Andrew Y.Z. Wang\thanks{Most of the work of this paper was done while Wang was
a Visiting Scholar in the Faculty of I.T., Monash University, Oct.\ 2015 -- Oct.\ 2016,
funded by the China Scholarship Council (CSC).}   \\
School of Mathematical Sciences   \\
University of Electronic Science and Technology of China   \\
Chengdu 611731   \\
P.R. China   \\
Email: \href{mailto:yzwang@uestc.edu.cn}{\texttt{yzwang@uestc.edu.cn}}
}
\date{16 May 2017}
\begin{document}
\maketitle

\begin{abstract}
A set $S\subseteq\{0,1\}^E$ of binary vectors, with positions indexed by $E$, is said to be
a \textit{powerful code} if, for all $X\subseteq E$, the number of vectors in $S$ that are zero
in the positions indexed by $X$ is a power of 2.  By treating binary vectors as characteristic
vectors of subsets of $E$, we say that a set $S\subseteq2^E$ of subsets of $E$ is a
\textit{powerful set} if the set of characteristic vectors of sets in $S$ is a powerful code.
Powerful sets (codes) include cocircuit spaces of binary matroids (equivalently,
linear codes over $\mathbb{F}_2$), but much more besides. Our motivation is that,
to each powerful set, there is an associated nonnegative-integer-valued rank function (by a construction
of Farr), although it does not in general satisfy all the matroid rank axioms.

In this paper we investigate the
combinatorial properties of powerful sets.  We prove fundamental results on special
elements (loops, coloops, frames, near-frames, and stars), their associated types of
single-element extensions, various ways of combining powerful sets to get new ones,
and constructions of nonlinear powerful sets.
We show that every powerful set is determined by its clutter
of minimal nonzero members.  Finally, we show that the number of powerful sets is doubly
exponential, and hence that almost all powerful sets are nonlinear.
\end{abstract}

%-------------------------------------------------------------

\section{Introduction}
\label{sec:intro}

Let $E$ be a finite set, called the \textit{ground set}, and let
$S\subseteq\{0,1\}^E$ be a set of binary vectors, with positions indexed by $E$.
A set $X\subseteq E$ of positions has the \emph{power-of-2 property (for $S$)} if
the number of vectors in $S$ that are zero on $X$
(i.e., in the positions indexed by $X$) is a power of 2.
We say $S$ is a \textit{powerful set}, or a \textit{powerful code},
if every $X\subseteq E$ has the power-of-2 property for $S$.
By treating binary vectors as characteristic
vectors of subsets of $E$, we also say that a set $S\subseteq2^E$ of subsets of $E$ is a
\textit{powerful set} if the set of characteristic vectors of sets in $S$ is a powerful set.
We move freely between subsets $X$ of $E$ and their characteristic vectors $\mathbf{x}$.
We prefer powerful \textit{set} terminology, but sometimes use powerful \textit{code}
terminology when commenting on connections with coding theory.

Unless stated otherwise, we use the ground set $E=[n]:=\{1,2,\ldots,n\}$.  We view $S$ as
a subset of the $n$-dimensional linear space $\mathbb{F}_2^n$, over the finite field
$\mathbb{F}_2$ consisting of all $01$-vectors of length $n$.

The \emph{order} of a powerful set $S$ is the size of its ground set, or equivalently, the
length of its vectors (when $S$ is viewed as a code).
The \emph{size} of $S$ is the cardinality of $S$.
The power-of-2 property for $X=E$ implies that the zero vector must be in $S$.
With $X=\emptyset$, we conclude that the size of $S$ is also a power of $2$.
The \emph{dimension} of $S$, written $\dim S$, is the nonnegative integer $d$ such that the size of $S$ is $2^d$.

Two powerful sets $S_1$ and $S_2$
are said to be \emph{isomorphic}, written $S_1\cong S_2$,
if there is a bijection between their
ground sets which induces a bijection between $S_1$ and $S_2$.

If $S$ is a finite-dimensional linear space over $\mathbb{F}_2$, then the vectors of $S$ that are 0 on $X$ form a subspace of $S$, thus the number of such vectors is a power of $2$.  Hence a linear space is always a powerful set.  From now on, we say a powerful set $S$ is \emph{linear} if it is a linear space, otherwise it is \emph{nonlinear}.  Up to isomorphism, there is a unique
smallest nonlinear powerful set, namely
\[
S=\{000, 011, 101, 111\}.
\]
Later we will see that almost all powerful sets are nonlinear.

For the sake of convenience, we often write a set $S\subseteq\mathbb{F}_2^n$ in the form
of a matrix whose rows are the elements of $S$.
For example, we can identify the above smallest nonlinear powerful set $S$ with the matrix
\begin{equation*}
\left(
\begin{array}{ccc}
0 &0 &0\\
0 &1 &1\\
1 &0 &1\\
1 &1 &1
\end{array}
\right).
\end{equation*}
We emphasise that, when $S$ is linear, this is not just a generator matrix for $S$; its
rows list \textit{all} members of $S$.

Our remarks above show that powerful sets generalise binary matroids, or equivalently,
binary linear codes.  Every binary matroid has a rank function $\rho:2^E\rightarrow\mathbb{N}\cup\{0\}$, defined on subsets of its ground set $E$, that satisfies the matroid rank
axioms.  Our original motivation for studying powerful sets was that they, too, have a
nonnegative-integer-valued ``rank-like'' function.  We elaborate on this now, before
setting the scene for the rest of this paper.

Let $f:2^E\rightarrow\{0,1\}$ be the indicator function of
a binary code $S\subseteq\mathbb{F}_2^E$, defined
for any $X\subseteq E$ by $f(X)=1$ or 0 according as the characteristic vector of $X$
does, or does not, belong to $S$.  The \textit{rank transform} $Q$, introduced in
\cite{farr93} (see also the exposition in \cite[\S3.6]{farr07a} and a closely related construction
due to Kung \cite{kung80}),
associates to any such $f$ the function $Qf$ defined on subsets of $E$ by
\begin{equation}
\label{eq:Qf}
Qf(X) = \log_2 \left(
\frac{\sum_{Y\subseteq E} f(Y)}{\sum_{Y\subseteq E\setminus X} f(Y)}
\right) .
\end{equation}
Observe that, when $Qf(X)$ is defined, it must be nonnegative.
(This follows from the fact that $f$ itself is nonnegative-valued.)
When it is defined, we call $Qf(X)$ the \textit{rank} of $X$, but bear in mind that this is
a loosening of that term since $Qf$ may not satisfy the matroid rank axioms.
For the special case when $S$ is linear, $Qf$ gives the usual rank function for the
binary matroid.  If $S$ is nonlinear, then $Qf$ may take irrational values
or be undefined for some arguments.
For $Qf(X)$ to be defined for all $X\subseteq E$, it is necessary and sufficient that
$f(\emptyset)=1$.  In particular, $Qf(X)$ is always defined if $S$ is a powerful set,
since in that case $\emptyset\in S$ so $f(\emptyset)=1$.
For $Qf$ to be integer valued, it is necessary and sufficient that $S$ be a powerful set.

Given that the functions $Qf$ extend rank functions, it is natural to investigate
what happens when they are used in place of rank functions.  This was done in
\cite{farr93,farr04}, where a theory of Tutte-Whitney polynomials is developed for
arbitrary functions $f:2^E\rightarrow\mathbb{C}$ (called \textit{binary functions}).
There, $Qf$ is used in place
of a matroid rank function to generalise the rank generating function
of Whitney \cite{whitney32} to arbitrary binary functions.
A surprising amount of Tutte-Whitney polynomial theory
extends to these objects, including duality, deletion-contraction relations, and
interesting partial evaluations.  But the ``polynomials'' themselves often have nonintegral
exponents.  It is therefore natural to focus on cases where the \emph{polynomials} are just
that, which means that $Qf$ is integer valued.  If, in addition, we ask that $f$ be
$\{0,1\}$-valued, so that it is indeed an indicator function and can be taken to represent
a subset of $2^E$, then we are led to the study of powerful sets.

If $S$ is a powerful set, we write $f_S$ for its indicator function, and $\rho_S$ for its
rank function, $\rho_S:=Qf_S$.

The definition of powerful codes is somewhat reminiscent of almost affine codes,
introduced in \cite{simonis-ashikhmin98}, although they are different in nature.
A $q$-ary code $S$ with index set $E$
is \textit{almost affine} if for all $X\subseteq E$ the cardinality of the code
$S_{E\backslash X} := \{(a_i)_{i\in E\backslash X}\mid (a_i)_{i\in E}\in S\}$ is a power of $q$.
The construction of $S_{E\backslash X}$ from $S$ is called \textit{puncturing}
with respect to $X$, or \textit{projection} onto $E\backslash X$.  We
simply discard all coordinates with positions in $X$, thereby shortening the vectors
to length $|E\setminus X|$.  This contrasts with powerful sets, where we do not remove
any coordinates, but simply require that the coordinates indexed by $X$ are zero.
When $q=2$, a binary code containing the zero vector is almost affine if and only if
it is linear \cite{simonis-ashikhmin98}, so binary almost affine codes give us nothing new,
and correspond to binary matroids.
See \cite{simonis-ashikhmin98,westerbaeck-etal2016} for further information
about almost affine codes and their connections with matroid theory.

In this paper we lay the foundations of the theory of powerful sets. We first
(in \S\ref{sec:reductions}) extend the contraction operation, for binary matroids,
to powerful sets.  Then, in \S\ref{sec:extns-special-elts}, we consider five types of
special elements: loops, coloops, frames, near-frames, and stars.  Of these, only loops
and coloops occur in binary matroids.  Each of the five has an associated type of
single-element extension operation, and we also generalise parallel extensions from binary
matroids to powerful sets. In \S\ref{sec:posn-wise-max}, we present a construction for some
nonlinear powerful sets, analogous to generating linear spaces from sets of vectors but using
positionwise maximum instead of positionwise addition in $\mathbb{F}_2$ (i.e., positionwise OR instead of
positionwise XOR). In \S\ref{sec:combining} we give three ways of combining powerful sets to
form new powerful sets.  Two of these have no real analogue for linear spaces.
Then in \S\ref{sec:generation} we show that every powerful set is determined by its clutter
of minimal nonzero members, by giving an algorithm to construct it from that clutter.
Finally, we consider enumeration of powerful sets in \S\ref{sec:enumeration}.  We report
the numbers of powerful sets (and, in particular, the numbers of nonlinear powerful sets)
of each order $\le6$.  The trend in this data is that nonlinear powerful sets quickly dominate,
and we confirm this trend mathematically.  We show that the number of loopless frameless
nonlinear powerful sets of order $n\ge5$ is
doubly exponential --- specifically, at least $2^{2^{(n-7)/3}}$ --- from which it follows that,
asymptotically, almost all powerful sets are nonlinear.

%--------------------------------------------------------------------------------------------------------------------

\section{Reductions}
\label{sec:reductions}

Let $S\subseteq2^E$ and $e\in E$.  Put
\[
S/e := \{ X\subseteq E\backslash\{e\} \mid X\in S \} .
\]
We say that $S/e$ is formed from $S$ by \textit{contraction} of $e$.
In terms of matrix representation, we remove column $e$ and also remove
all rows that have a 1 in the position indexed by $e$.
%\[
%S/e=\{(x_i)_{i\in E\backslash\{e\}}\mid (x_i)_{i\in E}\in S, x_e=0\}.
%\]

For example, consider the (nonlinear) powerful set $S = \{000, 011, 110, 111\}$, with
the usual ground set $\{1,2,3\}$.  Then
\[
\begin{array}{ccll}
S/1  & = &  \{00, 11\},  &  \hbox{with ground set $\{2,3\}$};   \\
S/2  & = &  \{00\},  &    \hbox{with ground set $\{1,3\}$};   \\
S/3  & = &  \{00, 11\},  &    \hbox{with ground set $\{1,2\}$}.
\end{array}
\]
So $S/1\cong S/3$.

\begin{thm}\label{contraction}
(a) If $S$ is powerful then $S/e$ is powerful.   \\
(b) If $S$ is linear then $S/e$ is linear.  (See, e.g., \cite[Theorem 9.3.1]{welsh76}.)
\qed
\end{thm}

The converses are not true, since (for example) adding a new all-0 column, indexed by $e$,
to $S/e$ (using the matrix representation viewpoint), then adding a row that is all-0 across
$E\backslash\{e\}$ but has 1 in position $e$, does not in general give another powerful set
(let alone a linear one).

The rank of $S/e$ is given by $\rho_{S/e}(X)=\rho_S(X\cup\{e\})-\rho_S(\{e\})$;
see \cite[\S4]{farr93}.

Another way of reducing a powerful set by a single element is to simply
delete the column indexed by $e$, without deleting any rows.
We call this \textit{deletion}, since it generalises deletion in binary matroids, and
denote the subset of $2^{E\backslash\{e\}}$ so formed by $S\backslash e$.
But, if it is applied to a nonlinear powerful set, it may leave duplicate rows in the reduced matrix,
giving a powerful multiset but not necessarily a powerful set.  The operation of \textit{puncturing}
with respect to $e$ consists of deletion of $e$ followed by removal of one member of each
pair of identical rows.  This ensures that we obtain a set rather than a multiset, and it yields a
linear powerful set if $S$ is linear (see, e.g., \cite{welsh76}), but it does
not necessarily produce a powerful set if $S$ is nonlinear.
Note also that the addition of a new column to a powerful set (i.e., the reverse of puncturing)
does not necessarily give a powerful set.

\section{Extensions and special elements}
\label{sec:extns-special-elts}

We now look at several ways to extend a powerful set by a single element.
A special role is played by five types of special elements.
The proofs are straightforward and most are omitted.

An element $e\in E$ that belongs to no set in $S$
(equivalently, it indexes a zero column in the matrix representation) is a \textit{loop},
and has rank 0.
The operation of adding a zero column to $T\subseteq\mathbb{F}_2^n$
is called \textit{loop extension}, and the resulting subset of $\mathbb{F}_2^{n+1}$
is denoted by $T+\circ$. Observe that, if $e$ is a loop of $S$,
then $S\backslash e=S/e$.

\begin{thm}\label{thm00}
If $e\in E$ is a loop of $S$ and $S/e$ is powerful then $S$ is powerful.
\qed
\end{thm}

Suppose that, writing $e$ as the last column and reordering rows if necessary,
$S\subseteq\mathbb{F}_2^n$ has a matrix of the form
\[
\left(
\begin{array}{cc}
T & \mathbf{0}   \\
T & \mathbf{1}
\end{array}
\right) ,
\]
where $T\subseteq\mathbb{F}_2^{n-1}$, and $\mathbf{0}$ and $\mathbf{1}$ are column vectors
whose length equals the size of $T$.  Then $e$ is a \textit{coloop} of $S$, and has rank 1.
The operation of forming $S$ from $T$ in this way is called \textit{coloop extension}.
We write $S=T+\circ^*$.

\begin{thm}\label{coloop}
If $e\in E$ is a coloop of $S$ and $S/e$ is powerful then $S$ is powerful.
\end{thm}

\pf For any $X\subseteq E$, if $e\not\in X$, then the number of vectors of $S$ that are 0 on $X$ is twice the number of vectors of $S/e$ that are 0 on $X$, thus being a power of $2$.

If $e\in X$, then the number of vectors of $S$ that are 0 on $X$ is the same as the number of vectors of $S/e$ that are 0 on $X\backslash\{e\}$, which is also a power of $2$.\qed

\begin{pro}
If $e\in E$ is a coloop of $S$ and $S/e$ is linear then $S$ is linear.
\end{pro}

\pf  For convenience, we write $e$ as the last column in the matrix representation of $S$.
Let $\mathbf{u}i$ and $\mathbf{v}j$ be any two vectors of $S$ where $\mathbf{u},\mathbf{v}\in S/e$ and $i,j\in\{0,1\}$.  It follows from the linearity of $S/e$ that $\mathbf{w}=\mathbf{u}+\mathbf{v}\in S$.  Thus we have $\mathbf{w}0\in S$ and $\mathbf{w}1\in S$.  Since $i+j\in\{0,1\}$, we can conclude that
\begin{align*}
\mathbf{u}i+\mathbf{v}j=\mathbf{w}(i+j)\in S.
\end{align*}
Thus $S$ is linear. \qed

\textsc{Remark.}  For a powerful set $S$, the zero row vector $\mathbf{0}$ belongs to $S$, thus $\mathbf{0}1\in S+\circ^*$.  Therefore, a coloop extension of a powerful set must have a vector of weight $1$.

\begin{con}\label{conjdup}
If $T$ is a powerful set with at least one vector of weight $1$, then $T$ is a coloop extension
of some powerful set $S$.
\end{con}

\textsc{Remark.} The conjecture is true for the linear case, since a singleton member of the
cocircuit space of a binary matroid must be a coloop.

Let $S$ be a powerful set, again with $e$ indexing the last column in its matrix, and now
with matrix of the form
\[
\left(
\begin{array}{cc}
\mathbf{0} & 0   \\
T\backslash\{\mathbf{0}\} & \mathbf{1}
\end{array}
\right) ,
\]
where $T$ is a powerful set, $\mathbf{0}$ is a row vector, and $\mathbf{1}$ is a column vector.  Note that $S\backslash e=T$.
Then $e$ is a \textit{frame} of $S$ (using terminology for an analogous concept
in \cite{welsh-whittle99}), and adjoining $e$ to $T$ is called \textit{framing} $T$ by $e$.
A frame has rank equal to $\dim S$.  We write $S=T+\Box$.

\begin{thm}\label{frame}
Let $S\subseteq\mathbb{F}_2^E$ have a frame $e\in E$.  Then $S\backslash e$ is powerful if and only if
$S$ is powerful.
\qed
\end{thm}

A powerful set can also be enlarged by an element that is almost, but not quite, a frame.

Suppose $S\subseteq\mathbb{F}_2^n$ and $\mathbf{v}\in S$ is a nonzero vector.
The set $S+\Box\backslash\mathbf{v}$ is formed by adding a new coordinate $0$ to the zero
vector $\mathbf{0}$ and $\mathbf{v}$, and a new coordinate $1$ to the remaining vectors of $S$.
The new element is called a \textit{near-frame} and has rank $\dim S-1$.

\begin{thm}
If $S$ is powerful then $S+\Box\backslash\mathbf{v}$ is powerful.
\qed
\end{thm}

If $T\subseteq\mathbb{F}_2^n$,
define $T+\star\subseteq\mathbb{F}_2^{n+1}$ by
\begin{align*}
T+\star=\{\mathbf{v}0\,|\,\mathbf{v}\in T\}\cup\{\mathbf{v}1\,|\,\mathbf{v}\notin T\}.
\end{align*}
We call the new element a \textit{star}.  If $T$ is powerful then the star has rank $n-\dim T$.

\begin{thm}\label{thmcon}
$T\subseteq\mathbb{F}_2^n$ is powerful if and only if $T+\star$ is powerful.
\end{thm}

\pf Let $\overline{T}$ be the set $\mathbb{F}_2^n\backslash T$.  Identify $S=T+\star$ with the following matrix
\begin{equation*}
S=\left(
\begin{array}{ccc|c}
{} &{} &{} &0\\
{} &T &{} &\vdots\\
{} &{} &{} &0\\\hline
{} &{} &{} &1\\
{} &\overline{T} &{} &\vdots\\
{} &{} &{} &1
\end{array}
\right).
\end{equation*}
For any $X\subseteq[n+1]$, if $n+1$ is not in $X$, then the number of rows that are 0
on $X$ must be a power of $2$. This is because the submatrix consisting of the first $n$ columns is the linear space $\mathbb{F}_2^n$. If $n+1\in X$, we only need to consider the submatrix
\begin{equation*}
S_T=\left(
\begin{array}{ccc|c}
{} &{} &{} &0\\
{} &T &{} &\vdots\\
{} &{} &{} &0
\end{array}
\right).
\end{equation*}
The rows of $S_T$ that are 0 on $X\backslash\{n+1\}$ are precisely those that are 0 on $X$. So the number of such rows is a power of 2 if and only if $T$ is a powerful set. Therefore $S$ is powerful if and only if $T$ is powerful. \qed

\begin{con}\label{conj}
Suppose that $S$ is a subset of $\mathbb{F}_2^n$ with $2^{n-1}$ elements, where $n\geq2$. If $S$ is a powerful set, then we can find a coordinate such that deleting this coordinate from all the elements of $S$ yields the set $\mathbb{F}_2^{n-1}$, i.e., all the new vectors are distinguishable.
\end{con}

\textsc{Remark.}  Conjecture \ref{conj} holds if $S$ is linear, since in that case we have
a binary matroid of rank $n-1$ on $n$ elements, which must have a circuit, and deleting
any element $e$ in the circuit gives a binary matroid $S\backslash e$ of rank $n-1$ on $n-1$
elements, whose cocircuit space is all of $\mathbb{F}_2^{n-1}$.
For the nonlinear case, Conjecture \ref{conj} holds for $n\le6$.

\textsc{Remark.} If we do not require that the size of $S$ is $2^{n-1}$, Conjecture \ref{conj} fails to hold. For example, let \[S=\{00000,00111,01011,01111,10101,10111,11010,11011\}.\]
It is easy to check that $S$ is a powerful set, but deleting any one coordinate will always yield two indistinguishable vectors of length $4$.

\textsc{Remark.} If a powerful set $S$ satisfies Conjecture \ref{conj}, it can always be constructed as $T+\star$ from a smaller powerful set $T$.  Suppose that deleting the last bit of each vector in $S$ gives all possible vectors of $\mathbb{F}_2^{n-1}$.  Collecting those vectors of $S$ whose last bit is $0$ and removing the last bit from each such vector yields the desired smaller powerful set.

If $S$ is a powerful set and $e\in E$, then the \textit{parallel extension}
of $S$, denoted by $S^{\shortparallel e}$, is formed by duplicating the column indexed by $e$ in the matrix representation of $S$.

\begin{thm}
Let $S\subseteq2^E$ and $e\in E$.  Then $S$ is powerful if and only if
its parallel extension $S^{\shortparallel e}$ is powerful.
\qed
\end{thm}

From binary matroid theory, we have

\begin{pro}
Let $S$ be a powerful set, then $S^{\shortparallel e}$ is linear if and only if $S$ is linear.
\qed
\end{pro}

%  Is this worth defining?
%If $S\subseteq2^E$ and $e\in E$, then the \textit{series extension} is constructed
%on the new ground set $E':=E\cup\{e'\}$, where $e'\not\in E$, as follows.
%\[
%S^{(\hbox{\scriptsize s}:e)} := \{\emptyset, \{e,e'\} \} \cup ((S\backslash\{\emptyset\})\times\{e,e'\}) .
%\]
%

%--------------------------------------------------------------------------------------------------------------------

\section{Position-wise max construction}
\label{sec:posn-wise-max}

Given any $S\subseteq\mathbb{F}_2^n$, elementary linear algebra gives us the
linear powerful set $\langle S\rangle$ consisting of all binary linear combinations of vectors
in $S$.  In this section we give another way to generate larger sets from $S$,
using a positionwise operation, which in this case will often give us nonlinear powerful sets.

A \emph{permutation matrix} is a square binary matrix that has exactly one entry of $1$ in each row and each column, and $0$s elsewhere.

Given any $S\subseteq\mathbb{F}_2^n$, we consider its matrix representation. If the matrix representation of $S$ contains a submatrix which is a permutation matrix of order $|S|$, then we say that $S$ is \emph{permutative}.

\textsc{Remark.} It is clear that a permutative set cannot contain the zero vector.

Define the \emph{disjunction} $\mathbf{u}\vee\mathbf{v}$ of two vectors
$\mathbf{u}=(u_i)_{i=1}^{n}$ and $\mathbf{v}=(v_i)_{i=1}^{n}$ in $\mathbb{F}_2^{n}$
by $\mathbf{u}\vee\mathbf{v}=(\mathrm{max}\{u_i,v_i\})_{i=1}^{n}$.

Suppose $S=\{\mathbf{u}_1,\mathbf{u}_2,\ldots,\mathbf{u}_m\}\subseteq\mathbb{F}_2^n$, define the \textit{disjunctive closure} of $S$ to be the set
\begin{align*}
\langle S\rangle_{\vee}=\{a_1\mathbf{u}_1\vee a_2\mathbf{u}_2\vee\cdots\vee a_m\mathbf{u}_m\,|\,a_i\in\mathbb{F}_2, 1\leq i\leq m\},
\end{align*}
where $a_i\mathbf{u}_i=\mathbf{u}_i$ if $a_i=1$, and $\mathbf{0}$ otherwise.
Note that the zero vector always belongs to $\langle S\rangle_{\vee}$.

\begin{thm}
If $m\leq n$ and $S=\{\mathbf{u}_1,\mathbf{u}_2,\ldots,\mathbf{u}_m\}\subseteq\mathbb{F}_2^n$ is a permutative set, then $\langle S\rangle_{\vee}$ is a powerful set of size $2^{m}$.
\end{thm}

\pf Since we are not concerned with order on $S$ or its ground set, we can assume that,
in the matrix representation
\begin{equation*}
S=\left(
\begin{array}{c}
\mathbf{u}_1   \\
\mathbf{u}_2   \\
\vdots         \\
\mathbf{u}_m
\end{array}
\right),
\end{equation*}
the first $m$ columns form the identity matrix $I_{m}$.

We first prove that any vector in $\langle S\rangle_{\vee}$ has a unique expression as $a_1\mathbf{u}_1\vee a_2\mathbf{u}_2\vee\cdots\vee a_m\mathbf{u}_m$, which shows that the size of $\langle S\rangle_{\vee}$ is $2^m$. Given a vector $\mathbf{v}=v_1v_2\cdots v_n\in \langle S\rangle_{\vee}$, we claim that
\begin{align*}
\mathbf{v}=v_1\mathbf{u}_1\vee v_2\mathbf{u}_2\vee\cdots\vee v_m\mathbf{u}_m.
\end{align*}
That is to say, $\mathbf{v}$ is completely determined by its first $m$ components. For $1\leq i\leq m$, $\mathbf{u}_i$ is the only vector in $S$ whose $i$th component is $1$. So if $v_i=1$, the coefficient of $\mathbf{u}_i$ must be $1$ otherwise the $i$th component of $\mathbf{v}$ will be $0$. Similarly, if $v_i=0$, the coefficient of $\mathbf{u}_i$ is $0$.

Next we show that $\langle S\rangle_{\vee}$ is a powerful set.  Given $X\subset[n]$, let $\mathbf{u}_{1,X},\mathbf{u}_{2,X},\ldots,\mathbf{u}_{r,X}\in S$ be all vectors that are 0 on $X$.  Then we claim that
\begin{align*}
\langle S\rangle_{\vee,X}:=\{a_1\mathbf{u}_{1,X}\vee a_2\mathbf{u}_{2,X}\vee\cdots\vee a_r\mathbf{u}_{r,X}\,|\,a_i\in\mathbb{F}_2,1\leq i\leq r\}
\end{align*}
contains all the vectors of $\langle S\rangle_{\vee}$ that are 0 on $X$.  It is clear that any vector $\mathbf{w}\in \langle S\rangle_{\vee,X}$ is 0 on $X$.  On the other hand, if
$\mathbf{w}\in \langle S\rangle_{\vee}$ is 0 on $X$, then in the unique expression
\begin{align*}
\mathbf{w}=w_1\mathbf{u}_1\vee w_2\mathbf{u}_2\vee\cdots\vee w_m\mathbf{u}_m,
\end{align*}
the coefficient of every $\mathbf{u}_i$ which is nonzero in some position in $X$ must be zero, otherwise $\mathbf{w}$ has a nonzero entry in some position in $X$.  Hence, the claim holds.  In addition, any two vectors of $\langle S\rangle_{\vee,X}$ are different, thus the size of $\langle S\rangle_{\vee,X}$ is $2^r$, a power of $2$.  If $X=[n]$, the zero vector is the only vector in $\langle S\rangle_{\vee}$ with all zero coordinates.  Therefore, $\langle S\rangle_{\vee}$ is a powerful set.\qed

\begin{exa}
\label{eg:perm}
Let $S=\{00011,01100,10101\}\subseteq\mathbb{F}_2^5$. Then the $1$st, $2$nd and $4$th columns of
\begin{equation*}
S=\left(
\begin{array}{ccccc}
0 &0 &0 &1 &1  \\
0 &1 &1 &0 &0  \\
1 &0 &1 &0 &1
\end{array}
\right)
\end{equation*}
comprise a permutation matrix of order $3$, thus $S$ is permutative.  We have
\begin{align*}
\langle S\rangle_{\vee}=\{00000,00011,01100,10101, 01111, 10111, 11101, 11111\}.
\end{align*}
It is straightforward to check that $\langle S\rangle_{\vee}$ is a powerful set.
\end{exa}

If $S$ is not permutative, then $\langle S\rangle_{\vee}$ is not necessarily powerful.  For example,
let $S=\{0111, 1011, 1101\}$, whose matrix representation has no unit vector columns
so $S$ is certainly not permutative.  Then $\langle S\rangle_{\vee}=\{0000, 0111, 1011, 1101, 1111\}$,
which has size 5, so is not powerful.

%--------------------------------------------------------------------------------------------------------------------

\section{Combining two powerful sets}
\label{sec:combining}

Basic set operations do not necessarily preserve the powerful property.  The
complement of a powerful set is never powerful (since it does not contain the zero vector),
and the union and intersection of powerful sets are not necessarily powerful.  (For example,
take the linear powerful set $\{000, 011, 101, 110\}$ and our smallest nonlinear powerful set $\{000, 011, 101, 111\}$.)

We now present three ways to combine two powerful sets which give another powerful set
(always, or under mild conditions).  Only the first corresponds to a binary matroid operation.

Let $Q\subseteq\mathbb{F}_2^{m}$ and $R\subseteq\mathbb{F}_2^{n}$.
The \textit{direct sum} of $Q$ and $R$ is defined by
\begin{align*}
Q\oplus R=\{\mathbf{u}\mathbf{v}\,|\, \mathbf{u}\in Q, \mathbf{v}\in R\}.
\end{align*}

\begin{thm}
$Q\oplus R$ is powerful if and only if $Q$ and $R$ are powerful.
\end{thm}
\pf If $X\subseteq[m]$ and $Y\subseteq\{m+1,\ldots,m+n\}$, then the number of vectors of $Q\oplus R$ that are zero on $X\cup Y$ is the number of vectors of $Q$ that are zero on $X$  times the number of vectors of $R$ that are zero on $Y$.  The result follows, paying particular
attention to the case $X=\emptyset$ and the case $Y=\emptyset$.\qed

Elementary linear algebra gives

\begin{pro}
The direct sum $Q\oplus R$ is linear if and only if $Q$ and $R$ are both linear.
\qed
\end{pro}

The direct sum generalises the direct sum of binary matroids and is a special
case of the product of disjoint binary functions \cite[p.\ 276]{farr93}.

We now come to our second way of combining powerful sets.

Write $\mathbf{0}_k$ and $\mathbf{1}_k$ for the row vector of $k$ $0$s and $k$ $1$s, respectively. Given two vectors $\mathbf{u}$ and $\mathbf{v}$, let $\mathbf{u}\mathbf{1}_n$ be the vector formed by appending $n$ $1$s to $\mathbf{u}$, and $\mathbf{1}_m\mathbf{v}$ be the vector formed by inserting $m$ $1$s before the start of $\mathbf{v}$, i.e., prepending $m$ $1$s to $\mathbf{v}$.

Let $Q\subseteq\mathbb{F}_2^{m}$ and $R\subseteq\mathbb{F}_2^{n}$ be powerful sets. Define the set $Q\#R\subseteq \mathbb{F}_2^{m+n}$ as follows
\[Q\#R=\{\mathbf{0}_{m+n}\}\cup\{\mathbf{u}\mathbf{1}_n\,|\,\mathbf{u}\in Q\backslash\{\mathbf{0}_m\}\}\cup\{\mathbf{1}_m\mathbf{v}\,|\, \mathbf{v}\in R\backslash\{\mathbf{0}_n\}\}\cup\{\mathbf{1}_{m+n}\}.\]
The construction of $Q\#R$ can be depicted as
\begin{equation*}
\left(
\begin{array}{c|c}
\mathbf{0}_m                  &\mathbf{0}_n \\[5pt]
Q\backslash\{\mathbf{0}_m\}   &\mathbf{1}_{(|Q|-1)\times n}\\[5pt]
\mathbf{1}_{(|R|-1)\times m}                  &R\backslash\{\mathbf{0}_n\} \\[5pt]
\mathbf{1}_m                  &\mathbf{1}_n
\end{array}
\right),
\end{equation*}
where $\mathbf{1}_{a\times b}$ is the all-one matrix with $a$ rows and $b$ columns.
\begin{exa}
If $Q=\{00\cdots0,11\cdots1\}\subseteq\mathbb{F}_2^m$ and $R=\{00\cdots0,11\cdots1\}\subseteq\mathbb{F}_2^n$, then
\[Q\#R=\{00\cdots0,11\cdots1\}\subseteq\mathbb{F}_2^{m+n},\]
which is also a powerful set.
\end{exa}

The result of combining powerful sets using \# is in general not powerful.
But there are many cases where it is,
and furthermore it can be used to construct nonlinear powerful sets.

\begin{thm}
Let $Q\subseteq\mathbb{F}_2^{m}$ and $R\subseteq\mathbb{F}_2^{n}$.
Then
$Q\#R$ is a powerful set if and only if $Q$ and $R$ are both powerful
and one of the following holds:
\begin{itemize}
\item[(a)]  one of $Q,R$ consists only of a zero vector and possibly an all-one vector, while the other includes an all-one vector; or
\item[(b)]  $|Q|=|R|$, and neither $Q$ nor $R$ contains an all-one vector.
\end{itemize}

Furthermore, if $Q\#R$ is powerful, then $Q\#R$ is nonlinear unless $Q$ and $R$ each consist just of a zero vector and possibly an all-one vector.
\end{thm}

\pf If $X\subseteq[m]$ is nonempty, then the vectors of $Q\#R$ that are $0$ on $X$ are precisely the vectors of $Q\backslash\{\mathbf{0}_m\}$ that are $0$ on $X$, each extended by $1$s at the end, together with $\mathbf{0}_{m+n}$. The number of these vectors is a power of $2$
if and only if $Q$ is a powerful set.

Similarly, if $Y\subseteq\{m+1,\ldots,m+n\}$ is nonempty, then the number of vectors of $Q\#R$ that are $0$ on $Y$ is a power of $2$ if and only if $R$ is a powerful set.

If $X\subseteq[m]$ and $Y\subseteq\{m+1,\ldots,m+n\}$, with each of $X$ and $Y$ being nonempty, then the only vector that is $0$ on $X\cup Y$ is $\mathbf{0}_{m+n}$, so the number is $2^0=1$.

Finally, the total number of vectors in $Q\#R$ (corresponding to the empty subset of positions) is
\[1+(|Q|-1)+(|R|-1)+1=|Q|+|R|,\]
provided $\mathbf{1}_m\not\in Q$ and $\mathbf{1}_n\not\in R$.
Under this condition,
if $Q$ and $R$ are powerful, then $|Q|=|R|$ if and only if
$|Q|+|R|$ is a power of $2$ if and only if $Q\#R$ is a powerful set.

Now suppose that $Q$ and $R$ are powerful, and either
$\mathbf{1}_m\in Q$ or $\mathbf{1}_n\in R$.
If just one of these holds then $|Q\#R|=|Q|+|R|-1$, which is not
a power of 2 unless exactly one of $|Q|,|R|$ is 1.  (They cannot both be 1, since one of $Q,R$
contains an all-one vector as well.)  In that case, the other is some power of 2 other than 1.
Suppose without loss of generality that $Q$ contains an all-one vector while $R$ contains only
a zero vector.  Then $Q\#R$ is equivalent to adding $n$ frames to $Q$.
If both $\mathbf{1}_m\in Q$ and $\mathbf{1}_n\in R$ then $|Q\#R|=|Q|+|R|-2$.  In that case,
one of $Q,R$ --- suppose $R$, without loss of generality ---
consists only of a zero vector and an all-one vector.  Again, we find that
$Q\#R$ is equivalent to adding $n$ frames to $Q$.
In any case, $Q\#R$ is powerful, by Theorem \ref{frame}.

We now consider nonlinearity.

If $Q$ and $R$ each consist just of a zero vector and possibly an
all-one vector, then $Q\#R$ consists just of the all-0 vector and the all-1 vector, so is
trivially linear.

Suppose then that (without loss of generality) $Q$ contains a vector $\mathbf{u}$ that is
nonzero and not all-ones.
We know that $\mathbf{u}\mathbf{1}_n$ is in $Q\#R$. It is clear that the last $n$ coordinates of $\mathbf{u}\mathbf{1}_n+\mathbf{1}_{m+n}$ are all $0$.
Since $\mathbf{u}\not=\mathbf{1}_m$ (as $\mathbf{1}_m\not\in Q$),
$\mathbf{u}\mathbf{1}_n+\mathbf{1}_{m+n}\neq\mathbf{0}_{m+n}$. But $\mathbf{0}_{m+n}$ is the unique vector in $Q\#R$ whose last $n$ coordinates are all $0$, which implies that $\mathbf{u}\mathbf{1}_n+\mathbf{1}_{m+n}\not\in Q\#R$. Therefore, $Q\#R$ is not a linear space.
\qed

It is interesting to consider the relationship between the rank functions
$\rho_Q, \rho_R, \rho_{Q\#R}$ of $Q,R,Q\#R$ respectively, when $Q\#R$ is powerful.
As for any powerful set, the
empty set has rank 0.  Now suppose $X, Y\not=\emptyset$, $X\subseteq[m]$ and $Y\subseteq\{m+1,\ldots,m+n\}$.  Then $\rho_{Q\#R}(X)=\rho_Q(X)+1$, $\rho_{Q\#R}(Y)=\rho_R(Y)+1$,
and $\rho_{Q\#R}(X\cup Y)=\dim Q+1=\dim R+1$.
In the light of this last observation, we call
$Q\#R$ the \textit{mutual framing} of $Q$ and $R$.

For $Q\subseteq\mathbb{F}_2^{n}$ and $R\subseteq\mathbb{F}_2^{n}$, define
\[
Q\bullet R:=\{\mathbf{v}00: \mathbf{v}\in Q\cap R\}\cup\{\mathbf{v}01:\mathbf{v}\in Q\backslash R\}\cup\{\mathbf{v}10:\mathbf{v}\in R\backslash Q\}\cup\{\mathbf{v}11:\mathbf{v}\not\in Q\cup R\}.
\]

\begin{thm}\label{thmgen}
$Q\bullet R$ is also powerful if and only if
$Q$, $R$ and $Q\cap R$ are all powerful.
\end{thm}
\pf For convenience, we identify $S=Q\bullet R$ with the matrix
\begin{equation*}
S=\left(
\begin{array}{c|c}
Q\cap R             &\mathbf{0}\mathbf{0} \\[8pt]
Q\backslash R       &\mathbf{0}\mathbf{1} \\[8pt]
R\backslash Q       &\mathbf{1}\mathbf{0} \\[8pt]
\overline{Q\cup R}  &\mathbf{1}\mathbf{1}
\end{array}
\right),
\end{equation*}
where $\overline{Q\cup R}$ is the complement of $Q\cup R$ in $\mathbb{F}_2^n$.

Consider any $X\subseteq[n+2]$.  We analyse whether it has the power-of-2 property
in the following four cases.

Case 1. If $n+1\not\in X$ and $n+2\not\in X$, then the number of rows which are $0$ on $X$ is a power of $2$ since the first $n$ columns of $S$ form the linear space $\mathbb{F}_2^n$.

Case 2. If $n+1\in X$ and $n+2\in X$, then the rows of $S$ that are 0 on $X$ are precisely those of $Q\cap R$ that are $0$ are $X\setminus\{n+1,n+2\}$, each extended by two $0$s at the end.
The number of these rows is a power of $2$ for all such $X$
if and only if $Q\cap R$ is a powerful set.

Case 3. If $n+1\in X$ and $n+2\not\in X$, we only need to consider the submatrix
\begin{equation*}
\left(
\begin{array}{c|c}
Q\cap R             &\mathbf{0} \\[8pt]
Q\backslash R       &\mathbf{0}
\end{array}
\right).
\end{equation*}
Since $Q=(Q\cap R)\cup(Q\backslash R)$, it follows that
$Q$ is a powerful set if and only if the number of rows that are $0$ on $X$ is a power of $2$
for all such $X$.

Case 4. If $n+1\not\in X$ and $n+2\in X$, the argument is similar to Case 3. \qed

\begin{exa}
Let $n=3$, and $Q=\{000, 001, 010, 011\}$ and $R=\{000,011,101,111\}$. It is easy to see that $Q$ and $R$ are powerful sets and $Q\cap R=\{000,011\}$ is also a powerful set. According to the construction in Theorem \ref{thmgen}, we have
\[Q\bullet R=\{00000, 01100, 00101, 01001, 10110, 11110, 10011, 11011\}.\]
It is straightforward to verify that $Q\bullet R$ is a powerful set.
\end{exa}

Using the cases of the above proof to analyse rank, we find that, for any $X\subseteq[n]$,
\begin{eqnarray*}
\rho_{Q\bullet R}(X)  & = &  |X|,   \\
\rho_{Q\bullet R}(X\cup\{n+1\})  & = &  n-\dim Q+\rho_Q(X),   \\
\rho_{Q\bullet R}(X\cup\{n+2\})  & = &  n-\dim R+\rho_R(X),   \\
\rho_{Q\bullet R}(X\cup\{n+1,n+2\})  & = &  n-\dim(Q\cap R)+\rho_{Q\cap R}(X).
\end{eqnarray*}

\textsc{Remark.} Theorem \ref{thmgen} can be extended further, using all possible three-bit extensions of vectors, with three powerful sets $P,Q,R$ with the appropriate intersections also having the power-of-2 property.  Then it could be extended to an arbitrary number $k$ of extra bits, with the same number of powerful sets with the required properties being combined.

%--------------------------------------------------------------------------------------------------------------------

%\section{Duality??? *****}
%
%Let  C  be a clutter with ground set  E.  (So  C  is a collection of subsets of  E  such that no member of  C  is a subset of any other member of  C.)
%
%The dual clutter  C*  (also called the blocker of  C)  is the set of all subsets  X  of  E  such that (a)  X  meets every member of  C  (i.e., for all  Y  in  C,  $X \cap Y$  is nonempty), and (b) no proper subset of  X  meets every member of  C  (i.e., if  $X' \subset X$  and  $X' \not= X$, then there exists  Y  in  C  such that  X  and  Y  are disjoint).  It is easy to see that the minimality condition, (b), means that  C*  is also a clutter.  This is all well known.
%
%Now, suppose  C  is a clutter that happens to be the set of minimal nonempty members of some powerful set.
%
%Is  C*  necessarily also the set of minimal nonempty members of some powerful set?
%
%If Yes, then we have a kind of duality for powerful sets.
%If No, then it would be good to have a counterexample.

%--------------------------------------------------------------------------------------------------------------------

\section{Generation}
\label{sec:generation}

Recall that a \textit{clutter} (also called a \textit{Sperner family}) is
an antichain in $2^E$ under the subset order.

If $S\subseteq2^E$ then $S_{\scriptsize\hbox{min}}$ denotes the set of
its minimal nonempty members, which is a clutter.

\begin{thm}
Every powerful set is determined by its minimal nonempty members.
\end{thm}

\pf Consider the following algorithm, which takes a clutter $S_0\subseteq2^E$ as input.
We will show that either it detects that there is no powerful set $S$ such that
$S_{\scriptsize\hbox{min}}=S_0$, and rejects $S_0$, or it computes an indicator function
$f:2^E\rightarrow\{0,1\}$ for a set $S=\supp f$ which is the unique powerful set
such that $S_{\scriptsize\hbox{min}}=S_0$ (where $\supp f := \{X\subseteq2^E\mid f(X)\not=0\}$).
\begin{tabbing}
\renewcommand{\baselinestretch}{3.0}
\hspace*{1cm}\=\ \ \ \ \ \=\ \ \ \ \ \=\ \ \ \ \ \=\hspace{1.8in}\=   \\
~1.  \>  Input:  $S_0$   \\
~2.  \>  $f(\emptyset) := 1$   \\
~3.  \>  For each $k=1,\ldots,n$   \\
     \>  \{   \\
~4.  \>  \>  For each $X\subseteq E$ such that $|X|=k$   \\
     \>  \>  \{   \\
~5.  \>  \>  \>  If  $X\in S_0$, then put $f(X) := 1$   \\
~6.  \>  \>  \>  else if  $\sum_{Y\subset X} f(Y) = 1$   \>  \>  //~   There is no $Y\subseteq X$ such that $Y\in S_0$.   \\[5pt]
~7.  \>  \>  \>  \>  $f(X) := 0$      \>  //~   This uses $X\not\in S_0$.   \\[5pt]
~8.  \>  \>  \>  else if  $\sum_{Y\subset X} f(Y) = 2$   \>  \>  //~   There is a unique $Y\subset X$ such that $Y\in S_0$.   \\[5pt]
~9.  \>  \>  \>  \>  $f(X) := 0$   \>  //~   To ensure $\sum_{Y\subseteq X} f(Y)$ is a power of 2.   \\[5pt]
10.  \>  \>  \>  else  \>  \>  //~   If we reach here, we know $\sum_{Y\subset X} f(Y)\ge3$.   \\
11.  \>  \>  \>  if $\sum_{Y\subset X} f(Y) = 2^i-1$ for some $i\ge2$   \\[5pt]
12.  \>  \>  \>  \>  $f(X) := 1$   \>  //~   To ensure $\sum_{Y\subseteq X} f(Y)$ is a power of 2.   \\[5pt]
13.  \>  \>  \>  else  \>  \>  //~   If we reach here, we know $\sum_{Y\subset X} f(Y)\ge4$.   \\
14.  \>  \>  \>  if $\sum_{Y\subset X} f(Y) = 2^i$ for some $i\ge2$   \\[5pt]
15.  \>  \>  \>  \>  $f(X) := 0$   \>  //~   To ensure $\sum_{Y\subseteq X} f(Y)$ is a power of 2.   \\[5pt]
16.  \>  \>  \>  else   \>  \>  //~   $\sum_{Y\subset X} f(Y)\not\in\{2^i-1,2^i\mid i\in\mathbb{N}\cup\{0\}\}$   \\[5pt]
17.  \>  \>  \>  \>  Reject $S_0$.  It cannot be $S_{\scriptsize\hbox{min}}$ for any powerful set $S$.   \\[5pt]
  \>  \>  \}   \\
  \>  \}   \\
18.  \>  Output $f$.   \\
19.  \>  Accept $S_0$.
\end{tabbing}

Suppose there exists a powerful set $S$ such that $S_{\scriptsize\hbox{min}}=S_0$.

We show by induction on $k$ that the above algorithm assigns, to all sets $X\subseteq E$
of size $k$, the value $f(X)=1$ if $X\in S$ and $f(X)=0$ otherwise.

Inductive basis:  for $k=0$, we have $X=\emptyset$, and the algorithm correctly assigns
$f(\emptyset)=1$ (in line 2) since $\emptyset\in S$.

Now let $k\ge1$ and
suppose the claim is true for all sizes $<k$, and let $X$ be any set of size $k$.

If $X\in S_0$, then the first condition of the cascaded if statement (line 5) is satisfied,
and the algorithm correctly sets $f(X)=1$.

Now suppose $X\not\in S_0$.

The order in which the algorithm visits the sets in $2^E$ ensures that it will visit all the proper
subsets $Y$ of $X$ before visiting $X$ itself.  When it reaches $X$,
it will have already assigned values $f(Y)$ to all $Y\subset X$.

By the inductive hypothesis, $\sum_{Y\subset X} f(Y)$ gives the number of proper subsets
of $X$ that belong to $S$.

So this sum equals 1 if and only if
no proper subset of $X$ is in $S$ except for $\emptyset$.
In this case, no proper subset of $X$
can be in $S_0$ either, by definition of $S$ and $S_0$.  So $X\not\in S$, else $X\in S_0$.
Now, in this case the algorithm takes the second option of the cascaded if statement (line 6)
and assigns $f(X)=0$ (in line 7), which is correct (in that $f$ is the indicator function of $S$ on this set
$X$).

It remains to consider cases where $\sum_{Y\subset X} f(Y)\ge2$, i.e.,
some nonempty proper subset of $X$ belongs to $S$.

The sum equals 2 if and only if there is exactly one nonempty proper subset of $X$ in $S$.
In this case, there are exactly two proper subsets of $X$ in $S$, which is already a power of 2,
so for $S$ to be powerful, we must have $X\not\in S$.
Here the algorithm takes the third option of the cascaded if statement (line 8), and correctly puts $f(X)=0$ (in line 9).

It remains to consider cases where $\sum_{Y\subset X} f(Y)\ge3$, i.e.,
the number of proper subsets of $X$ belonging to $S$ is at least 3.

If this quantity is one less than a power of 2, then in order for $S$ to be powerful, we must
have $X\in S$, and the algorithm takes the fourth option of the cascaded if statement
(lines 10--11) and correctly sets $f(X)=1$ (in line 12).

If this quantity equals a power of 2, then in order for $S$ to be powerful, we must
have $X\not\in S$, and the algorithm takes the fifth option of the cascaded if statement
(lines 13--14) and correctly sets $f(X)=0$ (in line 15).

Since we have assumed that $S$ is powerful, we know that the number of proper subsets of $X$
that belong to $S$ must be either a power of 2 or one less than a power of 2.  So the above
cases cover all possibilities, and the last option of the cascaded if statement (line 16) is never reached.
Therefore, we know that the algorithm always assigns the correct value $f(X)$ to $X$ so
that $f$ is the indicator function of $S$ on this set $X$.

Hence the claim is true, by induction.

Therefore, once the algorithm finishes, every $X\subseteq E$ will have been assigned a value
$f(X)$, and $f$ will be the indicator function of $S$.

Since the algorithm is deterministic, it finds (the indicator function of) the unique powerful set $S$
such that $S_{\scriptsize\hbox{min}}=S_0$.

If there is no powerful set $S$ such that $S_{\scriptsize\hbox{min}}=S_0$, then the algorithm
stops at a smallest set $X\subseteq E$
such that the sum $\sum_{Y\subset X} f(Y)\ge5$
and is neither a power of 2 nor one less than a power of 2.
It is impossible for any extension of
$f$ that includes $X$ in its domain to be the indicator function of a powerful set.
In this case, the algorithm takes the last option of the cascaded if
statement (line 16).  It does not assign a value to $f(X)$,
and it correctly rejects $S_0$ (in line 17).\qed

The clutter of minimal nonempty members of a powerful set plays a role
analogous to a basis of minimal vectors in a linear space.
Its members may be thought of as analogues, for powerful sets, of cutsets in graphs.

Some natural questions arise.
\begin{enumerate}
\item Can we characterise those clutters that consist of the minimal nonempty members of some powerful set?
\item What fraction of clutters come from powerful sets in this way?
\end{enumerate}

%--------------------------------------------------------------------------------------------------------------------

\section{Enumeration}
\label{sec:enumeration}

Let $p(n)$ be the number of isomorphism classes
of powerful sets of order $n$, and $\tilde{p}(n)$ be the number of isomorphism classes of nonlinear powerful sets of order $n$. By direct computation, with assistance from Peng Yang and Tingrui Yuan of UESTC, we have determined $p(n)$ and $\tilde{p}(n)$ for $n\le6$.
%\begin{table}[htb]
\begin{center}
{\renewcommand{\arraystretch}{1.3}
\begin{tabular}{|m{1.5cm}<{\centering}||m{1cm}<{\centering}|m{1cm}<{\centering}|m{1cm}<{\centering}|m{1cm}<{\centering}|m{1cm}<{\centering}|m{1cm}<{\centering}|}\hline
$n$      &$1$  &$2$  &$3$  &$4$ &$5$ &$6$ \\ \hline\hline
$p(n)$   &$2$  &$4$  &$9$  &$25$  &$102$ &$900$\\ \hline
$\tilde{p}(n)$  &$0$  &$0$  &$1$  &$9$   &$70$  &$832$ \\ \hline
\end{tabular}}%\caption{$p(n)$ and $\tilde{p}(n)$ for $1\leq n\leq 6$.}
\end{center}
%\end{table}

These numbers suggest that the number of powerful sets of order $n$ grows very rapidly
as $n$ increases, and that the proportion that are linear shrinks rapidly.

We now show that the number of isomorphism classes of nonlinear
powerful sets of order $n$ is doubly exponential in $n$, and
in fact this remains true if we restrict to size $2^{n-2}$.
It follows that almost all powerful sets are nonlinear.

To do this, we will use another way of combining powerful sets, based on operations
previously introduced.

Let $S_1,S_2\subseteq\mathbb{F}_2^{n}$ be powerful sets.  Define
$S_1\diamond S_2\subseteq\mathbb{F}_2^{n+3}$ by
\[
S_1\diamond S_2=(S_1+\circ)\bullet(S_2+\Box).
\]
This construction can be depicted as follows
\begin{equation*}
S_1\diamond S_2=\left(
\begin{array}{cccc}
  &00\cdots0 &\multicolumn{1}{|c}{0}  &\multicolumn{1}{|c}{00}\\[2pt] \hline
  &          &\multicolumn{1}{|c}{0}  &\multicolumn{1}{|c}{01}\\
  &S_1\backslash\{\mathbf{0}_n\} &\multicolumn{1}{|c}{\vdots} &\multicolumn{1}{|c}{\vdots}\\
  &          &\multicolumn{1}{|c}{0}  &\multicolumn{1}{|c}{01}\\ \hline
  &          &\multicolumn{1}{|c}{1}  &\multicolumn{1}{|c}{10}\\
  &S_2\backslash\{\mathbf{0}_n\} &\multicolumn{1}{|c}{\vdots} &\multicolumn{1}{|c}{\vdots}\\
  &          &\multicolumn{1}{|c}{1}  &\multicolumn{1}{|c}{10}\\ \hline
  &          &                        &\multicolumn{1}{|c}{11}\\
  &\mbox{the rest of $\mathbb{F}_2^{n+1}$}  &          &\multicolumn{1}{|c}{\vdots}\\
  &          &                        &\multicolumn{1}{|c}{11}\\
\end{array}
\right).
\end{equation*}

\begin{thm}
\label{thm:diamond-powerful}
If $S_1,S_2\subseteq\mathbb{F}_2^{n}$ are powerful sets, then $S_1\diamond S_2\subseteq\mathbb{F}_2^{n+3}$ is also a powerful set.
\end{thm}
\pf
Since $S_i$ (for $i=1,2$) is powerful, it follows from Theorem \ref{thm00} and Theorem \ref{frame} respectively that both $S_1+\circ$ and $S_2+\Box$ are powerful sets.  It is clear that $(S_1+\circ)\cap(S_2+\Box)=\{\mathbf{0}\}$, which is also a powerful set.  Now the desired result follows from Theorem \ref{thmgen}.
\qed

\begin{pro}\label{propn:diamond-loopless-frameless}
For any nontrivial powerful sets $S_1,S_2\subseteq2^E$, the set $S_1 \diamond S_2$ is loopless and frameless.
\end{pro}
%  YW:  The reason for nontriviality that if $S_i=\{\mathbf{0}\}$, then $S_1 \diamond S_2$ has a frame.

\pf It is clear from the construction that no loops or frames are created, regardless of $S_1$ and $S_2$.\qed

\begin{thm}
\label{thm:diamond-enum}
Let $\mathcal{S}$ be a set of nonisomorphic loopless frameless powerful sets of order $n$
and size $2^{n-2}$.  Then
\[
\mathcal{S}^{\diamond2} := \{ S_1 \diamond S_2 \mid S_1,S_2\in\mathcal{S} \}
\]
is a set of nonisomorphic loopless frameless powerful sets of order $n+3$ and size $2^{n+1}$.
If $n>3$, then every member of $\mathcal{S}^{\diamond2}$ is nonlinear.
\end{thm}

\pf Let $S_1,S_2\in\mathcal{S}$.  By Theorem \ref{thm:diamond-powerful}, $S_1 \diamond S_2$ is powerful.
By Proposition \ref{propn:diamond-loopless-frameless}, $S_1 \diamond S_2$ is loopless and frameless.

We now show that all the members of $\mathcal{S}^{\diamond2}$ are nonisomorphic.
Suppose, by way of contradiction, that there exist $S_1,S_2,S_1',S_2'\in\mathcal{S}$, with
either $S_1\not\cong S_1'$ or $S_2\not\cong S_2'$, such that
$S_1 \diamond S_2\cong S_1' \diamond S_2'$.
Let $f$ be an isomorphism from $S_1 \diamond S_2$ to $S_1' \diamond S_2'$.
Now, $f$ cannot map any element of $[n+1]$ to any element of $\{n+2,n+3\}$, since
the column $n+1+i$ has weight $2^{n+1}-|S_i|=2^{n+1}-2^{n-2}>2^n$ (for $i=1,2$),
while every column indexed by an $e\in[n+1]$ has weight $2^n$.

Let $i\in\{1,2\}$.  Since $S_i,S_i'$ are loopless and frameless,
$S_i\backslash\{\mathbf{0}\}$ and $S_i'\backslash\{\mathbf{0}\}$ each have no column
that is all-0 or all-1,
so they each have no column that looks like their portion of column $n+2$ or $n+3$.
There will therefore be only one other column that, in its rows corresponding to
$S_i\backslash\{\mathbf{0}\}$, matches column $n+2$ or $n+3$, and similarly for
$S_i'\backslash\{\mathbf{0}\}$: namely, column $n+1$.  Therefore $f(n+1)=n+1$.
We can now see that $f$ cannot mix $n+2$ from $n+3$ up, since the rows where
column $n+2$ is 0 are precisely the rows where column $n+1$ is 0, and the rows
where column $n+3$ is 0 are precisely the rows where column $n+1$ is 1.
So $f(n+2)=n+2$ and $f(n+3)=n+3$.

We have seen that $f$ maps $[n]$ to itself.  Also, for each $i=1,2$,
the mapping it induces on codewords of $S_1 \diamond S_2$
sends rows corresponding to $S_i$ to rows corresponding to $S_i'$ (else the last three
bits of the codewords do not match up).  Since (by assumption) $f$ is an isomorphism
from $S_1 \diamond S_2$ to $S_1' \diamond S_2'$, it must induce an isomorphism from
$S_1$ to $S_1'$ and from $S_2$ to $S_2'$.
Therefore $S_1\cong S_1'$ and $S_2\cong S_2'$.
This contradicts our assumption that $S_1\not\cong S_1'$ or $S_2\not\cong S_2'$.
(In fact, just one $S_i\cong S_i'$ is sufficient to get this contradiction.)

Therefore, all the members of $\mathcal{S}^{\diamond2}$ are nonisomorphic.

Since $n>3$, each $S_i$ has at least three nonzero members.  Let $\mathbf{x}$ and $\mathbf{y}$ be two nonzero members of $S_1$.  The corresponding vectors in $S_1 \diamond S_2$ have the same final three bits
(by construction), so their sum has last three bits all 0.  If $S_1 \diamond S_2$ is linear,
then this means that their sum is the $(n+3)$-bit zero vector,
since the only vector in $S_1 \diamond S_2$ with last two bits 0 is the zero vector.
This implies that $\mathbf{x}+\mathbf{y}=\mathbf{0}$.
This can only happen if $\mathbf{x}=\mathbf{y}$, which contradicts the
fact that they are distinct nonzero members of $S_1$.
Hence $S_1 \diamond S_2$ cannot be linear.
\qed

\begin{lem}
The number $q(n)$ of isomorphism classes of
loopless frameless nonlinear powerful sets of order $n\ge5$
and size $2^{n-2}$ satisfies $\log_2\log_2q(n)\ge(n-7)/3$.
\end{lem}

\pf We use induction on $n$.

For the base case, observe that there are at least two
nonisomorphic loopless frameless nonlinear powerful sets of order 5
and size $2^3$.  We saw one in Example \ref{eg:perm}, and another in the Remark
following Conjecture \ref{conj}.  It is therefore straightforward to construct two
nonisomorphic loopless frameless nonlinear powerful sets of any order $k$
and size $2^{k-2}$ (for example, using coloop extensions of the two of order 5 we have just
mentioned).
Therefore, for $k\in\{5,6,7\}$, we have
$q(k)\ge2$, so $\log_2\log_2q(k)\ge0\ge(k-7)/3$.

Now let $n\ge8$, and suppose that $\log_2\log_2q(k)\ge(k-7)/3$ for all $k$ such that
$5\le k<n$.  Let $\mathcal{S}$ be a set containing one representative of each isomorphism
class of loopless frameless nonlinear powerful sets
of order $n-3$ and size $2^{(n-3)-2}=2^{n-5}$.  By the inductive hypothesis,
$|\mathcal{S}|=:q(n-3)\ge2^{2^{(n-10)/3}}$.  By Theorem \ref{thm:diamond-enum},
$\mathcal{S}^{\diamond2}$ contains only loopless frameless nonlinear powerful sets of
order $n$ and size $2^{n-2}$, and they are all nonisomorphic.  We therefore have
\[
q(n)\ge|\mathcal{S}^{\diamond2}|=|\mathcal{S}|^2=q(n-3)^2\ge(2^{2^{(n-10)/3}})^2
=2^{2^{(n-7)/3}} .
\]

The result follows by induction.\qed

For an upper bound on $q(n)$, we can start with the number $2^{2^n}$
of all sets of subsets of $[n]$.  We saw in \S\ref{sec:generation}
that a powerful set is determined by its clutter of minimal nonempty members, so $q(n)$
is at most the number of inequivalent clutters of order $n$.
The number of clutters on $[n]$ is at least the number of sets of $\lfloor n/2\rfloor$-subsets of $[n]$, since any collection of distinct sets all of the same size is a clutter.
So the number of clutters is at least $2^{n\choose\lfloor n/2\rfloor}$.  Since each isomorphism
class of clutters has at most $n!$ members, the number of isomorphism classes of clutters
is at least $2^{n\choose\lfloor n/2\rfloor}/n!$.  This eventually
exceeds $2^{c^n}$ for any fixed $c<2$.
It follows that the number of inequivalent clutters does not give us a better upper bound
of the form $2^{c^n}$ than the na\"\i ve $2^{2^n}$.

The number of isomorphism classes of binary matroids on $n$ elements is well known to
satisfy the easy upper bound $2^{n^2}$.  It follows that, asymptotically, almost all powerful
sets are nonlinear.

%   **** comment on the data on this for small  n

%   **** define weight of vector

%--------------------------------------------------------------------------------------------------------------------
\section{Discussion}\label{sec:discussion}

We have laid some of the foundations of the theory of powerful sets, but there is much still
to be done.

One line of research is to consider aspects of binary matroid theory and determine
how far they extend to powerful sets.  Most of our work has been of this character,
including our Conjectures \ref{conjdup} and \ref{conj}.
In \S\ref{sec:generation} we proposed the problem of
characterising those clutters that are the set of minimal nonempty members of a powerful set,
which is analogous to characterising sets of circuits of binary matroids.
Research could also be done on Tutte-Whitney polynomials of powerful sets,
to determine what special properties they have beyond the general results of \cite{farr93,farr04}.

Another line of research is to examine the coding-theoretic properties of nonlinear powerful sets (viewed as powerful codes).  These are sufficiently general objects that many do not have
useful coding properties, but it is reasonable to expect that some classes of them
may be useful.

One could examine the relationship between linear codes over $\mathbb{Z}_4$ and the
binary codes obtained from them using the Gray map,
$0 \mapsto 00,  1 \mapsto 01,  2 \mapsto 11,  3 \mapsto 10$ (as suggested to us by
Peter Cameron).
This construction does not necessarily give a powerful set, as the following example shows.
On the left is a linear code over $\mathbb{Z}_4$ and on the right is the corresponding
binary code.
\[
\begin{array}{ccc}
000 & \mapsto & \underline{0}0\underline{0}0\underline{0}0   \\
013 & \mapsto & 000110   \\
022 & \mapsto & 001111   \\
031 & \mapsto & 001001   \\
101 & \mapsto & \underline{0}1\underline{0}0\underline{0}1   \\
110 & \mapsto & \underline{0}1\underline{0}1\underline{0}0   \\
123 & \mapsto & 011110   \\
132 & \mapsto & 011011   \\
202 & \mapsto & 110011   \\
211 & \mapsto & 110101   \\
220 & \mapsto & 111100   \\
233 & \mapsto & 111010   \\
303 & \mapsto & 100010   \\
312 & \mapsto & 100111   \\
321 & \mapsto & 101101   \\
330 & \mapsto & 101000
\end{array}
\]
For the binary code, the number of vectors that are 0 on $X=\{1,3,5\}$ is 3, not a power of 2.
(Note the underlined bits.)  So the binary code is not powerful.
It remains to determine which $\mathbb{Z}_4$-linear codes give nonlinear powerful codes,
and what properties they have.

Finally, we suggest the challenge of finding significantly stronger bounds on the number
(up to isomorphism) of powerful sets of order $n$, and determination of
\[\lim_{n\rightarrow\infty}(\log_2 p(n))^{1/n}.\]

\vspace{0.5cm}

\newpage
\noindent\textbf{Acknowledgements}

We thank Thomas Britz for helpful comments and drawing our attention to almost affine codes, Yongbin Li for some discussion, and Tingrui Yuan and Peng Yang for their assistance with the computations.  This work was
supported by the National Natural Science Foundation of China (No. 11401080).

%--------------------------------------------------------------------------------------------------------------------

\end{document}